\documentclass[11pt,a4paper]{article}
\usepackage{bbm,times,url,graphicx,anysize,amssymb,amsmath,amsthm}
\title{On Visibility Graphs of Point Sets in the Plane}
\author{\textsc{Florian Pfender}\thanks{Supported by the
  DFG Research Center {\sc Matheon} (FZT86)}\\
\normalsize TU Berlin, MA 6-2\\
\normalsize D-10623 Berlin, Germany\\
\normalsize \texttt{fpfender@math.tu-berlin.de}
\date{}}

\newcommand{\Z}{{\mathbbm Z}}

\newcommand{\R}{{\mathbbm R}}
\newcommand{\V}{{\mathcal V}}

\newtheorem{theorem}{Theorem}
\newtheorem{lemma}[theorem]{Lemma}

\newtheorem{corollary}[theorem]{Corollary}

\newtheorem{proposition}[theorem]{Proposition}

\begin{document}
\maketitle

\begin{abstract}
The visibility graph $\V(X)$ of a discrete point set $X\subset \R^2$ has vertex
set $X$ and an edge $xy$ for every two points $x,y\in X$ whenever there is no
other point in $X$ on the line segment between $x$ and $y$. We show
that for every graph $G$, there is a point set $X\in \R^2$, such that
the subgraph of $\V(X\cup \Z^2)$ induced by $X$ is isomorphic to
$G$. As a consequence, we show that there are visibility graphs of
arbitrary high chromatic number with clique number six settling a
question by K\'{a}ra, P\'{o}r and Wood.
\end{abstract}

\section{Introduction}

The concept of a visibility graph is widely studied in discrete
geometry. You start with a set of objects in some metric space, and
the visibility graph of this configuration contains the objects as
vertices, and two vertices are connected by an edge if the
corresponding objects can ``see'' each other, i.e., there is a
straight line not intersecting any other part of the configuration
from one object to the other. Often, there are extra restrictions on
the objects and on the direction of the lines of visibility.

Specific classes of visibility graphs which are well studied include
bar visibility graphs (see~\cite{H}), rectangle visibility graphs
(see~\cite{SW}) and visibility graphs of polygons (see~\cite{AK}). In this
paper we consider visibility graphs of point sets.

Let $X\subset \R^2$ be a a discrete point set in the plane. The {\em
visibility graph of $X$} is the graph $\V(X)$ with vertex set $X$ and 
edges $xy$ for every two points $x,y\in X$ whenever there is no 
other point in $X$ on the line segment between $x$ and $y$, i.e. when
the point $x$ is visible from the point $y$ and vice versa.

K\'{a}ra, P\'{o}r and Wood discuss these graphs \cite{KPW}, and make some observations
regarding the chromatic number $\chi(\V(X))$ and the clique number
$\omega(\V(X))$, the order of the
largest clique. In particular, they characterize all
visibility graphs with $\chi(\V(X))=2$ and $\chi(\V(X))=3$, and in
both cases, $\omega(\V(X))=\chi(\V(X))$. Similarly, they show the
following proposition.
\begin{proposition}\label{grid}
Let $\Z^2$ be the integer lattice in the plane, then 
$\omega(\V(\Z^2))=\chi(\V(\Z^2))=4$.
\end{proposition}
Note that $\V(\Z^2))$ is not perfect as it contains induced
$5$-cycles. Further, it is not true in general that
$\omega(\V(X))=\chi(\V(X))$---there are point sets with as few as nine
points with $\omega(\V(X))=4$ and $\chi(\V(X))=5$.

For general graphs, there are examples with $\chi(G)=k$ and
$\omega(G)=2$ for any $k$, one famous example is the sequence of
graphs $M_{k-2}$ by Mycielski~\cite{M}. No similar construction is
known for visibility graphs with bounded clique number. As their main result, K\'{a}ra
et al. construct a family of point sets with 
$\chi(\V(X))\ge (c_1\log~ \omega(\V(X_i)))^{c_2\log ~\omega(\V(X_i))}$ for
some constants $c_1$ and $c_2$ and with $\omega(\V(X_i))$ getting
arbitrarily large.
Our main result is the following theorem.
\begin{theorem}\label{main}
For every graph $G$, there is a set of points $X\subset \R^2$ such that
the subgraph of $\V(X\cup \Z^2)$ induced by $X$ is isomorphic to
$G$.
\end{theorem}
Let $G_k$ be a graph with $\chi(G_k)=k$ and
$\omega(G_k)=2$, and let $X_k$ be the corresponding set given by Theorem~\ref{main}.
Let $Y_k\subset X_k\cup \Z^2$ be the subset of points contained in the
convex hull of $X_k$. Then $\chi(\V(Y_k))\ge \chi(G_k)=k$ and
$\omega(\V(Y_k))\le \omega(G_k)+\omega(\V(\Z^2))=6$, so we get the
following corollary settling the 
question from above raised by K\'{a}ra et al.
\begin{corollary}
For every $k$, there is a finite point set $Y\subset \R^2$, such that $\chi(\V(Y))\ge k$ and
$\omega(\V(Y))=6$.
\end{corollary}

\section{Proof of the Theorem}

Let $G$ be a graph with vertex set $V(G)=\{ 1,2,\ldots,n\}$ and edge
set $E(G)$. We will show the
following lemma in the Section~\ref{sec3}.
\begin{lemma}\label{lem}
For $M$ large enough, there is a set of prime numbers $\{ p_{ij}:1\le i<j\le n\}$ with the
following properties:
\begin{enumerate}
\item $2^M<p_{ij}<2^{M+1}$.
\item For $1\le k\le n$, let $P_k=2^{M+k+n_k}\prod_{i=1}^{k-1} p_{ik}\prod_{j=k+1}^n p_{kj}$,
  with $0\le n_k\le k+n-1$ chosen such that $\log_2 P_k=nM+2k$. 
Then $p_{k\ell}$ is the only number in $\{ p_{ij}:1\le i<j\le n\}$
  which divides $P_\ell-P_k$ for $1\le k<\ell\le n$.
\end{enumerate}
\end{lemma}
\noindent
From this, we can construct the set of points $X$ in
Theorem~\ref{main}:
$$X=\{ x_i:1\le i\le n\}\subset \R^2,
\mbox{ with }
x_i=\left( 2^{-nM}P_i~,~i~\frac{\prod\limits_{k<j}(P_j-P_k)}{\prod\limits_{kj\in
      E(G)}p_{kj}}\right) .$$
Before we prove the lemma, we will show that this point set has the
properties stated in the theorem.
For $1\le i<\ell\le n$, let $m_{i\ell}$ be the slope of the line through
$x_i$ and $x_\ell$. Then
$$m_{i\ell}=\frac{\ell-i}{P_\ell-P_i}\cdot\frac{2^{nM}\prod\limits_{k<j}(P_j-P_k)}{\prod\limits_{kj
    \in E(G)}p_{kj}}.$$
There are no three colinear points in $X$, as 
$$2^{nM+2i+1}\le P_{i+1}-P_{i}<2^{nM+2i+3},$$
thus $m_{i(i+1)}>m_{(i+1)(i+2)}$, and therefore $m_{i\ell}>m_{ik}$ for
$i<\ell<k$. Thus, $\V(X)$ is complete, and it remains to show that there
is an integer point on the line segment between $x_i$ and $x_\ell$ if and
only if $i\ell\notin E(G)$. To establish this goal, we look at the
intersections of the line through $x_i$ and $x_\ell$ ($i<\ell$) with the gridline
with constant $x$-coordinate $s\in \Z$, where $1\le
2^{-nM}P_i<s<2^{-nM}P_\ell\le 2^{2n}$. For this 
intersection point $z_{i\ell}^s=(s,y_{i\ell}^s)$ we have
\begin{multline*}
y_{i\ell}^s=i~\frac{\prod\limits_{k<j}(P_j-P_k)}{\prod\limits_{kj\in E(G)}p_{kj}}+(s-
2^{-nM}P_i)m_{i\ell}\\
= \underbrace{i~\frac{\prod\limits_{k<j}(P_j-P_k)}{\prod\limits_{kj\in E(G)}p_{kj}}}_{(1)}+
\underbrace{s~\frac{\ell-i}{P_\ell-P_i}\cdot\frac{2^{nM}
    \prod\limits_{k<j}(P_j-P_k)}{\prod\limits_{kj\in E(G)}p_{kj}}}_{(2)}+ 
\underbrace{P_i~\frac{\ell-i}{P_\ell-P_i}\cdot\frac{\prod\limits_{k<j}(P_j-P_k)}{\prod\limits_{kj\in
      E(G)}p_{kj}}}_{(3)}.
\end{multline*}
The expression (1) is an integer since $p_{kj}$ divides
$P_j-P_k$. By the same argument , (3) is an integer---just note
further that $p_{i\ell}$ divides $P_i$. It remains the analysis of
(2). 

If $i\ell\notin E(G)$, then (2) is an integer. Therefore,
$z_{i\ell}^s\in \Z^2$, and
$x_ix_\ell\notin E(\V(X\cup \Z^2))$.
If $i\ell\in E(G)$, observe that $p_{i\ell}>2^M>\max\{ \ell-i,s\}$, so $p_{i\ell}$
does not divide $s$ or $\ell-i$. Clearly, $p_{i\ell}$ does not divide
$2^{nM}$, and by Lemma~\ref{lem}, it does not divide any of the
$P_j-P_k$ other than $P_\ell-P_i$. Thus, (2) is not an integer,
$z_{i\ell}^s\notin \Z^2$, and 
$x_ix_\ell\in E(\V(X\cup \Z^2))$, proving Theorem~\ref{main}.\qed

\section{Proof of Lemma~\ref{lem}}\label{sec3}

By an inequality of Finsler~\cite{F}, there are
more than
$2^M/(3(M+1)\ln 2)>2n^3$ prime numbers in the interval from $2^M$ to $2^{M+1}$.

We will pick the $p_{ij}$ sequentially in the order
$p_{12}, p_{13},\ldots, p_{1n}, p_{23},\ldots, p_{(n-1)n}$.
Assume that we have picked numbers up to but not including $p_{ij}$
according to the Lemma, and we want to pick $p_{ij}$.
If $j\ne n$, we pick any prime number in the interval from $2^M$ to
$2^{M+1}$ that was not selected before (condition (*)), such that $p_{ij}$ does not
divide $P_k-P_\ell$ for all $1\le \ell<k<i$ (condition (**)). There were less than
${n\choose 2}$ primes selected before, and each $P_k-P_\ell$ has at most $n$
prime divisors greater than $2^M$, thus at most ${n\choose 2}+n{n \choose
  2}<n^3$ of the choices are blocked, and so this choice is possible.

If $j=n$, we have to fulfill (*) and (**), with the
added condition (***) that no $p_{k\ell}$ may divide $P_i-P_r$ for $\{ k,\ell\}\ne
\{ i ,r\}$. Pick $p_{ij}$ according to (*) and (**), and assume that
$p_{k\ell}$ divides $P_i-P_r$ for some $\{ k,\ell\}\ne\{ i ,r\}$.
We have $k\ne i$ as all $p_{i\ell}$ divide $P_i$, otherwise
$p_{i\ell}$ also divides $P_r$ and thus $r=\ell$, a
contradiction. Similarly, $\ell\ne i$.

Pick another number $p_{ij}'$ according to (*) and (**). If
$p_{k\ell}$ divides $P'_i-P_r$, then $p_{k\ell}$ divides
$P'_i-P_i=(p'_{ij}-p_{ij})P_i/p_{ij}$, and thus $p_{k\ell}$ divides
$p'_{ij}-p_{ij}$. But this is impossible since
$|p'_{ij}-p_{ij}|<2^M<p_{k\ell}$.
Therefore, each $p_{k\ell}$ can block at most one choice for $p_{ij}$,
so in total at most ${n\choose 2}$ further choices are blocked by
(***), and we can always find a number $p_{ij}$ with (*), (**) and
(***).
This concludes the proof of the lemma.
\qed

\section{Further Questions}

We have shown that there are visibility graphs with $\chi(\V(X))\ge k$ and
$\omega(\V(X))=6$ for every $k$. For all visibility graphs with
$\omega(\V(X))\le 3$, we know that $\chi(\V(X))=\omega(\V(X))$. The
only cases left to consider are $\omega(\V(X))=4$ and
$\omega(\V(X))=5$. A similar technique of combining a visibility graph
with $\omega(\V(X))=3$ with a graph $G$ with $\omega(G)=2$ and large
chromatic number will not work, since the visibility graphs with
$\omega(\V(X))=3$ are too simple (all but at most two of their
vertices are collinear unless $\V(X)$ is a special graph on six
vertices). It would be no surprise to us if the chromatic number of
visibility graphs with $\omega(\V(X))=5$ is bounded.

Finally, one could look for smaller point sets with $\chi(\V(X))\ge k$ and
$\omega(\V(X))=6$, as our sets tend to be very large.

\bibliographystyle{amsplain}

\end{document}